\documentclass[11pt]{article}
\usepackage{fullpage}

\newtheorem{theorem}{Theorem}[section]
\newtheorem{lemma}[theorem]{Lemma}

\newtheorem{corollary}[theorem]{Corollary}

\newenvironment{proof}[1][Proof]{\begin{trivlist}
\item[\hskip \labelsep {\bfseries #1}]}{\end{trivlist}}
\newenvironment{definition}[1][Definition]{\begin{trivlist}
\item[\hskip \labelsep {\bfseries #1}]}{\end{trivlist}}

\newenvironment{remark}[1][Remark]{\begin{trivlist}
\item[\hskip \labelsep {\bfseries #1}]}{\end{trivlist}}

\newcommand{\qed}{\nobreak \ifvmode \relax \else
      \ifdim\lastskip<1.5em \hskip-\lastskip
      \hskip1.5em plus0em minus0.5em \fi \nobreak
      \vrule height0.75em width0.5em depth0.25em\fi}

\usepackage{amsfonts}
\usepackage{verbatim}

\def \cali {{\cal I}}

\def \bprf {\begin{proof} {\bf :} }
\def \eprf {\end{proof}}
\def \beqn {\begin{eqnarray}}
\def \eeqn {\end{eqnarray}}
\def \bthm {\begin {theorem}}
\def \ethm {\end {theorem}}
\def \bcorr {\begin{corollary}}
\def \ecorr {\end{corollary}}
\def \brem {\begin{remark}}
\def \erem {\end{remark}}

\begin{document}
\title{New lower bounds for the independence number of
 sparse graphs and hypergraphs}
\author{Kunal Dutta\thanks{The Institute of Mathematical Sciences, Taramani, Chennai - 600113, India, email: kdutta@imsc.res.in},\, Dhruv Mubayi\thanks{
Department of Mathematics, Statistics, and Computer Science,
University of Illinois, Chicago, IL-60607, USA
email: mubayi@math.uic.edu.  Research supported in part by NSF grant DMS 0969092},\, C.R. Subramanian\thanks{The Institute of Mathematical Sciences, Taramani, Chennai - 600113, India,  email: crs@imsc.res.in}}

\maketitle
\date{}

\abstract{
We  obtain new lower bounds for the independence number of
$K_r$-free graphs and linear $k$-uniform hypergraphs in terms of the
degree sequence. This answers some old questions raised by Caro and Tuza
\cite{CT91}.  Our proof technique is an extension of a method
of Caro and Wei \cite{CA79, WE79}, and we also give a new short proof
of the main result of \cite{CT91} using this approach.
As byproducts, we also obtain some non-trivial identities involving
binomial coefficients.

\section{Introduction } 
For $k \ge 2$, a $k$-uniform  hypergraph $H$ is a pair $(V(H),E(H))$ where
$E \subseteq {V(H) \choose k}$.  A set $I\subset V(H)$ is an independent set of $H$ if
$e \not\subseteq I$ for every $e \in E(H)$, or
equivalently, ${I \choose k} \cap E(H) = \emptyset$.
The independence number  of $H$, denoted
by $\alpha(H)$, is the maximum size of an independent set in $H$.
For $u \in V(H)$, its degree in $H$, denoted by $d_H(u)$, is defined to be
$|\{e \in E(H) : u \in e\}|$ (we omit the subscript if it is obvious from context).
{\bf Throughout} this paper, we use $t$ to denote $k-1$ except in some
places where it stands for some real value (the correct meaning can be
easily inferred from the context).
Also, we use the term graph whenever $k$ happens to be $2$.
    A $k$-uniform hypergraph is \emph{linear} if it has no 2-cycles where
a 2-cycle is a set of 2 hyperedges containing at most $2t$ vertices.
The dual of the above definition  says that a
linear hypergraph is one in which every pair of vertices is contained in at most one
hyperedge.

\paragraph{}
In \cite{TU41}, Tur\'{a}n proved a theorem giving a tight bound on the
maximum number of edges that a $K_r$-free graph can have, which has
since become the cornerstone theorem of extremal graph theory.
Tur\'{a}n's theorem, when applied to the
complement $\overline{G}$ of a graph $G$, yields a lower
bound  $\alpha(G) \ge \frac{n}{d+1}$ where $d$ denotes the
average degree in $G$ of its vertices.

\paragraph{}
Caro \cite{CA79} and Wei \cite{WE79} independently proved that
$\alpha(G)\ge \sum_{v} \frac{1}{d(v)+1}$ which is at least $\frac{n}{d+1}$.
The probabilistic proof of their result
later appeared in \cite{ASBook}.
\footnote{According to R. Bopanna
\cite{FOBlog}, the probabilistic argument in \cite{ASBook} was
obtained by him,
although it is possible that it was known earlier.}
One natural extension of Tur\'an's theorem to $k$-uniform hypergraphs
$H$ is the bound $\alpha(H)> c_k\frac{n}{d^{1/t}}$, and this was shown
via an easy probabilistic argument by Spencer~\cite{Spencer}.
Caro and Tuza \cite{CT91} improved this bound for irregular $k$-uniform
hypergraphs by proving that
\begin{equation} \label{ct}
\alpha(H) \ge \sum_{v \in V(H)} \frac{1}{{d(v) + 1/t \choose d(v)}}.
\end{equation}
Indeed, an easy consequence of (\ref{ct}) is the following result.
\begin{theorem} {\bf (Caro-Tuza \cite{CT91})} \label{caroweihypgr}
For every $k \geq 3$, there exists  $d_k>0$
such that every $k$-uniform
hypergraph $H$ has
$$\alpha(H) \geq d_k\sum_{v\in V(H)}\frac{1}{(d(v)+1)^{1/t}}.$$
\end{theorem}
As a corollary, one infers the bound of Spencer above.
Later, Thiele \cite{Th99} provided a lower bound on the
independence number of non-uniform hypergraphs, based on the
degree rank (a generalization of degree sequence).

\paragraph{}
In this paper, we prove new lower bounds for the independence number of locally sparse graphs and linear $k$-uniform hypergraphs.  The starting point of our approach  is the probabilistic proof of Boppana-Caro-Wei.  This approach, together with some additional simple ideas,  quickly yields a new short proof of Theorem \ref{caroweihypgr} (see Section 2 for the detailed proof).

\subsection{$K_r$-free graphs}
\parindent=0pt

For certain classes of sparse graphs, improvements of the Caro-Wei bound (in terms of average degree $d$) are known.
Ajtai, Koml\'{o}s and Szemer\'{e}di \cite{AKS81}
proved a lower bound of $\Omega\left(\frac{n\log d}{d}\right)$ for
the independence number of
triangle-free graphs. An elegant and simpler proof was later given by
Shearer \cite{Shearer83}, who also improved the constant involved.
Later Shearer \cite{Shearer95} also proved a bound of
$\Omega\left(\frac{n\log d}{d\log \log d}\right)$ for
$K_r$-free graphs when $r>3$.

\paragraph{}
Caro and Tuza \cite{CT91} raised the following question in their
1991 paper :

\medskip

\noindent
$(i)$ Can the lower bounds of Ajtai \emph{et al} \cite{AKS81}
and Shearer (\cite{Shearer83}, \cite{Shearer95}) be generalized in
terms of degree sequences?
\medskip

We answer this question via the following two theorems.

\begin{theorem} \label{trianglefreethm}
For every  $\epsilon \in [0,1)$ there exists $c>0$ such that the following holds:
Every triangle-free graph $G$ with average degree $D$ has independence number at least
$$c(\log D)\sum_{v\in V(G)} \frac{1}{\max{\{D^{\epsilon},d(v)\}}}.$$
\end{theorem}

\begin{theorem} \label{krfreethm}
For every  $\epsilon \in [0,1)$ and $r \ge 4$, there exists $c>0$ such that the following holds:
Every $K_r$-free graph $G$ with average degree $D$  has independence number at least
$$c\frac{\log D}{\log \log D}\sum_{v\in V(G)} \frac{1}{\max{\{D^{\epsilon},d(v)\}}}.$$
\end{theorem}

\subsection{Linear Hypergraphs}
As mentioned earlier, a lower bound of $\Omega\left(n/d^{1/t}\right)$   for an $n$ vertex $k$-uniform hypergraph with average degree $d$ can be inferred from Theorem \ref{caroweihypgr}.
Caro and Tuza~\cite{CT91} also raised the following
question:
  \medskip

 $(ii)$ How can one extend the lower bounds of Ajtai \emph{et al} \cite{AKS81}
and Shearer (\cite{Shearer83}, \cite{Shearer95}) to hypergraphs?
\smallskip

As it turns out, such extensions were known for the class of
linear $k$-uniform hypergraphs.  Indeed, the lower bound
\begin{equation} \label{dlr} \alpha(H) = \Omega\left(n\left(\frac{\log d}{d}\right)^{1/t}\right),\end{equation} where $H$ is a linear $k$-uniform hypergraph with average degree $d$
was proved by Duke-Lefmann-R\"odl \cite{DLR95}, using the results of \cite{AKPSS82}. Our final result generalizes (\ref{dlr}) in terms of the degree sequence
of the hypergraph.

\begin{theorem} \label{kunithm}
For every $k \geq 3$ and $\epsilon \in [0,1)$, there exists $c>0$ such that
the following holds: Every linear $k$-uniform hypergraph $H$ with average degree $D$ has independence
number at least
$$c(\log D)^{1/t}\sum_{v\in V(H)} \frac{1}{\max{\{D^{\epsilon/t},(d(v))^{1/t}\}}}.$$
\end{theorem}

We also describe an infinite family of $k$-uniform linear
hypergraphs to illustrate that the ratio between the bounds of
Theorem \ref{kunithm}
and (\ref{dlr}) can be unbounded in terms of the number
of vertices.

\paragraph{}
The remainder of this paper is organized as follows. In Section 2, we give a new short proof of Theorem \ref{caroweihypgr}.
In Section \ref{sec:linazerobds}, we apply the analysis in
Section 2 to the special case of linear hypergraphs, and obtain
a ``warm-up" result - Theorem \ref{azerothm}, which will
be helpful in proving the main technical result, Theorem
\ref{genrlathm}, proved in Section \ref{sec:lingenrlabds}.
The expression obtained in Theorem \ref{genrlathm}  plays a crucial role in the proofs of Theorems \ref{trianglefreethm}, \ref{krfreethm} and \ref{kunithm}; these are provided in Section 5. In Section 6,
we give infinite families of $K_r$-free graphs
and $k$-uniform linear hypergraphs which illustrate that the
bounds in Theorems \ref{trianglefreethm}, \ref{krfreethm} and \ref{kunithm}
can be bigger than the corresponding bounds in
\cite{AKS81,AKPSS82,DLR95,Shearer83,Shearer95} by arbitrarily large
multiplicative factors. Finally, in section \ref{sec:concln},
we state several combinatorial identities which follow as simple corollaries of Theorem \ref{genrlathm}.

\section{A new proof of Theorem \ref{caroweihypgr}}
In this section we obtain a new short proof of Theorem \ref{caroweihypgr}.
First we obtain the following theorem which is later used to prove
Theorem \ref{caroweihypgr}.
\begin{theorem} \label{hgraphineqthm}
For every $k \geq 2$, there exists a constant $c=c_k$ such that any
$k$-uniform hypergraph $H$ on $n$ vertices and $m \geq 1$ hyperedges satisfies
$$ \sum_{J\subset V(H)} \frac{1}{{n\choose {|J|}}} > c\frac{n}{m^{1/k}}
\;\;\; \ldots \ldots \;\;\; (A) $$
where we sum over all independent sets $J$.
\end{theorem}
\begin{proof}
Let $t_k(n,m)$ denote the LHS of $(A)$.
Consider any edge $e \in E(H)$. $e$ can belong to at
most ${{n-k}\choose{j-k}}$ non-independent sets of size $j$. Since there
are $m$ edges there are at most $m{{n-k}\choose{j-k}}$ sets of size $j$ that
are not independent. Thus, at least ${n\choose j} - m{{n-k}\choose{j-k}}$
sets of size $j$ are independent.
Hence we have
\begin{eqnarray*}
 t_k(n,m) & \geq &
\sum_{j=1}^n \left(1-m\frac{{{n-k}\choose{j-k}}}{{{n}\choose{j}}}\right)
\;\; = \;\;  \sum_{j=1}^n \left(1-m\frac{(j)_k}{(n)_k}\right)  \\
&>& \sum_{j=1}^{\lfloor n/(2m)^{1/k} \rfloor}
\left(1-m\frac{j^k}{n^k}\right)
 \;\;  \geq  \;\;  \sum_{j=1}^{\lfloor n/(2m)^{1/k}\rfloor}
\left(1-m\frac{1}{2m}\right)  \\
&\geq & \frac{1}{2} \left \lfloor \frac{n}{(2m)^{1/k}}\right \rfloor
\;\;  \geq \;\;  c_k \frac{n}{m^{1/k}}
\end{eqnarray*}
for some suitably chosen $c_k$ which is close to $2^{-(k+1)/k}$.
\hfill $\Box$
\end{proof}
Let $H=(V,E)$ be a $k$-uniform hypergraph. For $k \geq 3$ and
for $u \in V$ with $d_H(u) \geq 1$,
the link graph associated with $u$ in $H$ is the
$t$-uniform hypergraph $L_u = (U,F)$ where
$U := \{v \not= u : \exists e \in E : \{u,v\} \subseteq e \}$ and
$F = \{e \setminus u : u \in e \in E\}$.  Let $\cali(H)$ denote the collection of independent sets of $H$.
\bigskip

{\bf Proof of Theorem \ref{caroweihypgr}.}
Let $H=(V,E)$ be an arbitrary $k$-uniform hypergraph.
Choose uniformly at random a total ordering $<$ on $V$.
Define an edge $e \in E$ to be {\em backward} for a vertex $v \in e$ if
$u < v$ for every $u \in e \setminus \{v\}$.
Define a random subset $I$ to be the set of those vertices $v$ such
that no edge $e$ incident at $v$ is backward for $v$ with respect to
$<$. Clearly, $I$ is independent in $H$. We have
$E[|I|] = \sum_{v} Pr(v \in I)$.  If $d_v=0$, then
$v \in I$ with probability $1$. Hence, we assume that $d(v) \geq 1$.
From the definition of $I$, it
follows that $v \in I$ if and only if for every $e$ incident at $v$,
$e \setminus \{v\} \not\subseteq S_v \:= \{u \in V(L_v) : u < v\}$.
In other words, $S_v$ is an independent set in $L_v$.
Let $l_v= |V(L_v)|$.
Then
\begin{eqnarray*}
Pr[v \in I] &=& \sum_{J \in \cali(L_v)}
 \frac{|J|!(l_v-|J|)!}{(l_v+1)!} \;\; =  \;\;
\frac{1}{l_v+1}\sum_{J \in \cali(L_v)}\frac{1}{{{l_v}\choose{|J|}}}
\end{eqnarray*}
Applying Theorem \ref{hgraphineqthm} to the $t$-uniform link graph
$L_v$ (with $c=c_{k-1}$), we get
\begin{eqnarray*}
    Pr[v\in I] &\geq& \frac{c}{l_v+1}\left(\frac{l_v}{d(v)^{1/(k-1)}}\right)
    \;\; \geq \;\; \frac{cl_v}{l_v+1}\left(\frac{1}{(d(v)+1)^{1/(k-1)}}\right).
\end{eqnarray*}
Since
$l_v\geq k-1$, we get $Pr[v\in I] \geq ((k-1)c/k)\frac{1}{(d(v)+1)^{1/(k-1)}}$.
By choosing $d_k=(k-1)c/k$, we get the lower bound of the theorem.\hfill
$\Box$

\section{Linearity : Probability of having no backward edges}
\label{sec:linazerobds}

In this section, we state and prove a warm-up result on the probability
of having no backward edges incident at a vertex for a randomly chosen
linear ordering (Theorem \ref{azerothm} below). The problem is the same as in the previous section, only,
now the hypergraph under consideration is assumed to be linear and we
get an explicit closed-form expression for this probability.
This result will be helpful for the proof of the main technical theorem,
given in the next section.
In order to state the lower bound, we need the following definition (of
fractional binomial coefficients) from \cite{GrKnPa94}.

\begin{definition}
For $t > 0$, $a \geq 0$, $d \in \mathbb{N}$
$${{d+1/t}\choose a} := \frac{(td+1)(t(d-1)+1)...(t(d-a+1)+1)}{a!t^a}$$
\end{definition}

\begin{theorem} \label{azerothm}
Let $H$ be a linear $k$-uniform hypergraph and let $v$ be an arbitrary
vertex having degree $d$. For a uniformly chosen total ordering $<$ on $V$,
the probability $P_v(0)$ that $v$ has no backward edge incident at it, is
given by $$P_v(0) = \frac{1}{{{d+1/t}\choose d}}$$
\end{theorem}

{\bf Remark.} It is interesting to note that the above expression when
summed over all vertices, is the same bound which Caro
and Tuza obtain in \cite{CT91} (using very different methods), although
their bound holds for independent sets in \emph{general} $k$-uniform
hypergraphs.

\paragraph{} We prove the theorem using the well-known Principle of
Inclusion and Exclusion (PIE). First we state an identity involving
binomial coefficients.

\begin{lemma} \label{azerolemma}
Given non-negative integers $d$ and $t$,
\begin{eqnarray*}
\sum_{r=0}^d (-1)^r{d\choose r}\frac{1}{tr+1} &=& \frac{1}{{{d+1/t}\choose d}}
\label{binident1}
\end{eqnarray*}
\end{lemma}
This identity is already known (see \cite{GrKnPa94}, Equation 5.41).
However, we give an alternate proof (using hypergeometric series) in the Appendix.

\begin{proof}[Proof of Theorem \ref{azerothm}]
Firstly, observe that since $H$ is linear, the number of
vertices that are neighbors of $v$ is exactly $(k-1)d = td$. Next, notice
that since the random ordering is uniformly chosen, only the relative
arrangement of these $td$ neighbors and the vertex $v$, i.e.
$td+1$ vertices in all, will
determine the required probability. Hence the total number of orderings
under consideration is $(td+1)!$.

\paragraph{}
Label the hyperedges incident at $v$ with $1,...,d$ arbitrarily. For a
permutation $\pi$, we say that $\pi$ has the property
$T_{\geq S}$ if
the edges with labels in $S$, $S \subseteq [d]$ are backward. Also, say
$\pi$ has the property $T_{=S}$ if the edges with labels
in $S$ are backward and no other edges are backward.
For a set $S$ of hyperedges incident at $v$, let
$N(T_{\geq S})$ denote the number of
orderings having the property $T_{\geq S}$, that is, the number of
permutations such that the hyperedges in $S$ will \emph{all} be backward
edges. $N(T_{= S})$ is similarly defined. $N(T_{\geq S})$ is determined
as follows : \\
Suppose $S$ has $r$ hyperedges incident at $v$.
For a \emph{fixed} arrangement of the vertices belonging to
edges in $S$, the number of
permutations of the remaining vertices is $(td+1)!/(tr+1)!$.
In each allowed permutation, the vertex $v$ must occur only after the
vertices of $S$
(i.e. the rightmost position). However the remaining $tr$ vertices can be
arranged among themselves in $(tr)!$ ways. Thus we have
$$N(T_{\geq S})= (td+1)!\frac{(tr)!}{(tr+1)!} = \frac{(td+1)!}{(tr+1)}.$$

\paragraph{}
Clearly, if a permutation has the property $T_{\geq S}$, it has the property
$T_{= S'}$ for some $S'\supseteq S$.
Hence for every $S \subset [d]$,
 $$N(T_{\geq S}) = \sum_{S' \supseteq S} N(T_{= S'}). $$
Therefore, by PIE (see \cite{STBook}, Chapter 2),
$$N(T_{=\emptyset}) = \sum_{S} (-1)^{|S|} N(T_{\geq S})$$

$$\sum_{|S|=r} N(T_{\geq S}) \;\; = \;\; {d \choose r} N(T_{\geq [r]}) \;\; =\;\;  {d\choose r} \frac{(td+1)!}{tr+1}$$
Hence we get the required probability to be
\begin{eqnarray*}
P_v(0) &=& \left(\sum_{r=0}^d {d \choose r}(-1)^r \frac{(td+1)!}{tr+1}\right)\times \frac{1}
             {(td+1)!}\\
    &=& \sum_{r=0}^d {d \choose r}(-1)^r \frac{1}{tr+1}
\end{eqnarray*}
By Lemma \ref{azerolemma},
$$P_v(0) = \frac{1}{{{d+1/t}\choose d}}$$
and this completes the proof. \hfill $\Box$
\end{proof}

\section{Linearity : Probability of having few backward edges}
\label{sec:lingenrlabds}

Now, we consider the more general case when at most $A-1$ backward
edges are allowed. In this section, we get an exact expression for
the corresponding probabiity. This estimate plays an important role
later in getting new and improved lower bounds on $\alpha(H)$ for
locally sparse graphs and linear hypergraphs.
Our goal in this section is to prove the following result.

\begin{theorem} \label{genrlathm}
For a $k$-uniform linear hypergraph $H$, a vertex $v$ having degree $d$,
a uniformly chosen permutation $\pi$ induces at most $A-1$ backward edges with
probability $P_v(A-1)$ given by
$$P_v(A-1) =  \left\{ \begin{array}{ll}
         1 & \mbox{if $d \leq A-1$};\\
         \frac{tA}{tA+1}\frac{{d\choose A}}{{{d+1/t}\choose d-A}} &
         \mbox{if $d \geq A$}.\end{array} \right.$$
\end{theorem}

\vspace{2mm}

\begin{corollary} \label{genrlacor}
As $d \to \infty$, the asymptotic expression for the probability $P_v(A-1)$
is given by
$$P_v(A-1) \;\;\sim\;\; \frac{1}{1+(1/(tA))}\left(\frac{A}{d}\right)^{1/t} \;\;=\;\;
                   \Omega(\left(A/d\right)^{1/t})$$
\end{corollary}

\begin{proof}
The asymptotics are w.r.t. $d \to \infty$, $d \geq A$.
The expression for having at most $A-1$ backward edges is
\begin{eqnarray*}
P_v(A-1) &=& \frac{1}{1+(tA)^{-1}}\frac{d(d-1)...(A+1)}{(d-A)!}\frac{(d-A)!}{(d+1/t)(d-1+1/t)...(A+1+1/t)} \\
&=& \frac{1}{1+(tA)^{-1}}\frac{1}{(1+1/td)(1+(t(d-1))^{-1})...(1+(t(A+1))^{-1})}
\end{eqnarray*}
Now, for $0< x$, we have $(1+x)^{-1} > e^{-x}$. So we get
\begin{eqnarray*}
P_v(A-1) &>& (1+(tA)^{-1})^{-1}e^{(-1/t)\sum_{r=1}^{d-A}(1/(A+r))} \\
         &=& (1+(tA)^{-1})^{-1}e^{(-1/t)[\sum_{r=1}^{d}(1/r)-\sum_{r=1}^A (1/r)]} \\
         &=& (1+(tA)^{-1})^{-1}e^{(-1/t)[\ln d- \ln A]+O((d-A)/(tdA))} \\
         &=& (1+(tA)^{-1})^{-1}e^{(-1/t)\ln (d/A)-O((d-A)/(tdA))} \\
         &=& (1+(tA)^{-1})^{-1}(A/d)^{1/t}\Omega(1)  \\
         &=& \Omega((A/d)^{1/t})
\end{eqnarray*}
The above expression therefore becomes $\Omega((A/d)^{1/t})$.
\hfill $\Box$
\end{proof}

\vspace{2mm}
The version of PIE used most commonly deals with $N(T_{= \emptyset})$,
i.e. the number of elements in the set of interest - in this case,
permutations of $[td+1]$ which do not have \emph{any} of the properties
under consideration (in this case, backward edges with respect to $v$).
However we need something slightly different - an expression for the
number of permutations which have \emph{at least} $A$ backward edges.
Clearly, the remaining permutations are those which have \emph{at most} $A-1$
backward edges.
\paragraph{}
Therefore, we use a slightly modified version of PIE, which is
stated below in Theorem \ref{mpiethm}.
This form is well-known (see e.g. \cite{STBook}, Chapter 2,
Exercise 1), although it seems to be used less frequently. For the sake
of completeness, we provide a simple proof.
First we state two identities involving binomial coefficients that we will prove in the Appendix.
\begin{lemma} \label{mpielemma}
For $a,b$ nonnegative integers,
$$\sum_{i=0}^b (-1)^i {{a+b}\choose{a+i}}{{a+i-1}\choose i} \;\; = \;\; 1$$
\end{lemma}

\begin{lemma} \label{genrlaexprlemma}
Given non-negative integers $d,A$, $d\geq A$ and a positive integer $t$,
$$\sum_{r=0}^{d-A} (-1)^r{d\choose {r+A}}{{A+r-1}\choose r}\frac{1}{t(r+A)+1} \;\;
=\;\; 1-\left(\frac{At}{tA+1}\right)\frac{{d\choose A}}{{{d+1/t}\choose d-A}}$$
\end{lemma}

We now present the generalized PIE and its well-known proof.

\begin{theorem} \label{mpiethm}
Let $S$ be an $n$-set and $E_1,E_2,...E_d$ not necessarily distinct subsets
of $S$. For any subset $M$ of $[d]$, define $N(M)$ to be the number of
elements of $S$ in $\cap_{i\in M} E_i$ and for $0\leq j \leq d$, define
$N_j := \sum_{|M|=j} N(M)$. Then the number $N_{\geq a}$ of elements of $S$ in
\emph{at least} $a$, $0\leq a\leq d$ of the sets $E_i$, $1\leq i\leq d$, is
$$ N_{\geq a} = \sum_{i=0}^{d-a} (-1)^i {{a+i-1}\choose i}N_{i+a} \;\;\;\; ...  \;\;\;\; (MPIE)$$
\end{theorem}

\begin{proof}
Take an element $e \in S$.
\begin{enumerate}
\item[(i)] Suppose $e$ is in no intersection of at least $a$ $E_i$'s.
Then $e$ does not contribute to any of the summands in the RHS of the
expression (MPIE), and hence, its net contribution to the RHS is zero.
\item[(ii)] Suppose $e$ belongs to exactly $a+j$ of the $E_i$'s, $0\leq j
\leq d-a$. Then its contribution to the RHS of (MPIE) is
$$\sum_{l=0}^j (-1)^l {{a+j}\choose {a+l}}{{a+l-1}\choose l}$$
and by Lemma \ref{mpielemma} this is equal to 1.\hfill $\Box$
\end{enumerate}

\end{proof}

\paragraph{}
\begin{proof}[Proof of Theorem \ref{genrlathm}]
  If $d \leq A-1$, then $P_v(A-1)=1$ obviously. The proof is similar to the proof of Theorem \ref{azerothm}, except that in
place of the PIE, we use Theorem \ref{mpiethm}. The set under consideration
is the set of permutations of $[td+1]$, the subsets $E_i$ correspond to
the permutations for which the $i$-th edge is backward. It is easy to see that
$N(M)= N(T_{\geq M})$ under the notation used in Theorem \ref{azerothm}
and hence $ N(M) = \frac{(td+1)!}{t|M|+1}$. Therefore we have
$N_j = {d\choose j}\frac{(td+1)!}{tj+1}$ as before.
Hence the expression for the probability $Q_v(A)$ that \emph{at least} $A$
edges \emph{are} backward under a uniformly random permutation $\pi$, becomes:
\begin{eqnarray*}
Q_v(A) &=& \sum_{i=0}^{d-A} (-1)^i{d\choose{i+A}}{{A+i-1}\choose i}
           \frac{1}{t(i+A)+1}.
\end{eqnarray*}
By Lemma \ref{genrlaexprlemma} the RHS of the above expression is
$$Q_v(A) \;\; = \;\; 1-\left(\frac{1}{1+(tA)^{-1}}\right)
    \frac{{d\choose A}}{{{d+1/t}\choose d-A}}.$$
Hence the probability of having at most $A-1$ backward edges is given by
$$P_v(A-1) \;\;=\;\; \frac{1}{1+(tA)^{-1}}
                     \frac{{d\choose A}}{{{d+1/t}\choose d-A}}$$
and the proof is complete. \hfill $\Box$\end{proof}

\section{Lower bounds for linear hypergraphs and
         $K_r$-free graphs} \label{sec:applcn}

In this section we prove Theorems  \ref{trianglefreethm}, \ref{krfreethm}, and \ref{kunithm}. These follow by a simple application of Corollary
\ref{genrlacor}.
Since the proofs follow the same outline, we prove them simultaneously,
highlighting only the differences as and when they occur.

\begin{proof} [Proofs of Theorems  \ref{trianglefreethm}, \ref{krfreethm} and \ref{kunithm}.]
Consider a uniformly chosen random permutation of the vertices of the
graph/hypergraph under
consideration. Let $D$ be the average degree of the graph or hypergraph and $A=D^{\epsilon}$. Let $I$ be the set of those vertices each having at most
$A-1$ backward edges incident on it. Clearly, the expected size of $I$ is

$$E[|I|] \;\;\geq \;\;c\sum_{v\in V} \left(\frac{A}{\max{\{A,d(v)\}}}\right)^{1/t}
         \;\;= \;\;cA^{1/t}\sum_{v\in V} \left(\frac{1}{\max{\{A,d(v)\}}}\right)
                   ^{1/t} $$

for some constant $c=c(k,\epsilon)$. 
(For a graph, $k=2$ and hence $t=1$). Also, by construction, the average
degree of the sub(hyper)graph induced by $I$ is at most $k(A-1)$.
Therefore, there exists an
independent set $I'$ of size at least as follows
\begin{enumerate}
\item[(i)] Case $t=1$, graph is $K_3$-free: By \cite{Shearer83}, $\alpha(G)$ is
at least
$$\Omega\left(\log (2(A-1))\frac{|I|}{2(A-1)}\right) \;\;=\;\; \Omega\left(\log D\sum_{v\in V} \frac{1}{\max{\{A,d(v)\}}}\right)$$
\item[(ii)] Case $t=1$, graph is $K_r$-free $(r>3)$: By \cite{Shearer95},
$\alpha(G)$ is at least
$$\Omega\left(\frac{\log (2(A-1))}{\log\log (2(A-1))}\frac{|I|}{2(A-1)}\right) \;\; = \;\;
                                    \Omega\left(\frac{\log D}{\log\log D}\sum_{v\in V}
                                    \frac{1}{\max{\{A,d(v)\}}}\right)$$
\item[(iii)] Case $t>1$, hypergraph is linear: By \cite{DLR95}, $\alpha(H)$
is at least
$$\Omega\left((\log k(A-1))^{1/t}\frac{|I|}{(k(A-1))^{1/t}}\right) \;\; = \;\;
                                    \Omega\left((\log D)^{1/t}\sum_{v\in V}
                                    \frac{1}{(\max{\{A,d(v)\}})^{1/t}}\right)$$
\end{enumerate}
The above three cases prove Theorems \ref{trianglefreethm}, \ref{krfreethm}
and \ref{kunithm} respectively.

\paragraph{}
{\bf Note}: An inspection of the proofs above show why we need $\epsilon$ to
be a fixed constant. It is because all three expressions above essentially have
$\log A$ i.e. $\epsilon\log D$ in the numerator. So, if $\epsilon=o(1)$, then
$\log A = o(\log D)$, and we would get asymptotically weaker results.
\hfill $\Box$
\end{proof}

\section{Construction comparing average degree vs.
degree sequence based bounds} \label{sec:exmpl}

A degree sequence-based bound obviously reduces to
a bound based on average degree, when the (hyper)graph
is regular.
However, the convexity of the function
$x^{-1/t}$, $x\geq 1$ and $t \in \mathbb{N}$, shows
that the bounds in Theorems \ref{trianglefreethm}, \ref{krfreethm}
and \ref{kunithm} are better
than the corresponding average degree-based bounds
proved in \cite{AKPSS82}, \cite{Shearer83} and
\cite{Shearer95} respectively provided the minimum degree is at least
$A$, although it is not
clear \emph{a priori} if the improvement can become significantly
larger. Also, at least half the vertices will have
degree at most $2D$, so even in the general
case (no restriction on the minimum degree) our bounds are no worse
than the average degree based bounds (ignoring the constant factors).
However, they can be much larger than the latter bounds.
We now give infinite families of $K_r$-free graphs and linear
$k$-uniform hypergraphs which show that
\begin{enumerate}
\item[(i)] The bounds given by Theorem \ref{trianglefreethm},
\ref{krfreethm} can be better than the bounds in \cite{AKS81, Shearer83,
Shearer95} respectively by a multiplicative factor of $\log (|V(G)|)$.
\item[(ii)] The bound in Theorem \ref{kunithm} can be better
than the bound in \cite{AKPSS82} by a multiplicative factor of
$\left((\log |V(H)|)/(\log \log |V(H)|)\right)^{(1-\epsilon)/t}$,
where $\epsilon$ is the constant mentioned in Theorem \ref{kunithm}.
\end{enumerate}

\begin{enumerate}
\item[Case (i)] Take a set of $n$ disjoint graphs, $K_{1,1}$,
$K_{2,2}$, $K_{4,4}$, ... , $K_{2^{n-1},2^{n-1}}$. The total
number of vertices is $2^{n+1}-2$, whereas the average degree is
$d_{av}=(2^n+1)/3$. Hence, the average degree based bound
gives $\Theta(|V(G)|\log d_{av}/d_{av}) = \Theta(\log d_{av})$.
Denote by $l$ the maximum $j$ such that $2^j \leq A= d_{av}^{\epsilon}$.
It follows that $A < 2^{l+1}$. Theorem \ref{trianglefreethm} gives
\begin{eqnarray*}
c\log d_{av}\sum_{v\in V} \frac{1}{\max\{d(v),A\}} & = &
c\log d_{av}\left[\sum_{j=0}^{l}\frac{2.2^j}{A} +
\sum_{j=l+1}^{n}\frac{2.2^j}{2^j}\right] \\
&=& c(\log d_{av})\left[\Theta(1) + \Theta(n)\right] \\
&=& c(\log d_{av})\Theta(\log(|V(G)|)
\end{eqnarray*}
The same example works for Theorem \ref{krfreethm} also, since
triangle-free graphs are obviously $K_r$-free, for $r\geq 3$.

\item[Case (ii)]
Fix some $m=m(n)=w \cdot k^{2^n}$, $w=\omega(1)$.
For each $i\in \{0,\ldots,n-1\}$, first create a $k$-uniform
$2^i$-regular linear hypergraph as follows: Take the
vertex set as $[k]^{2^i}$, and let each hyperedge consist of
the $k$ vertices which have all but one co-ordinate
fixed. Call this hypergraph an $i$-\emph{unit}. Now for each
$i$, create a component by taking
$\lceil\frac{m}{k^{2^i}}\rceil$ disjoint unions of
$i$-units. 
Take the disjoint unions of $n$ such components, one for every
$i \in \{0,\ldots,n-1\}$, to get the hypergraph
$H=H(n)=(V,E)$.
Now, the total number of vertices in the $i$-th component is
$m_i= m(1+o(1))$, and hence $|V|= nm(1+o(1))$. Also, the
average degree is $d_{av} \sim (2^n-1)/n \sim 2^n/n$. Let
$l$ denote the greatest integer $j$ such that $2^j \leq
(d_{av})^{\epsilon}\sim
2^{\epsilon n}/n^{\epsilon}$.
Therefore the average degree based bounds in
\cite{AKPSS82, DLR95} give a lower bound of

$$\alpha(H) = \Omega(mn^{1+1/t}(\log d_{av})^{1/t}/2^{n/t})\;\;\ldots\;\; (A)$$

On the other hand, the bound in Theorem \ref{kunithm} gives
\begin{eqnarray*}
\alpha(H) &=& \Omega\left((\log d_{av})^{1/t}\left[\sum_{j=0}^{l}
    \frac{mn^{\epsilon /t}}{2^{\epsilon n/t}}+\sum_{j=l+1}^{n-1}
    \frac{m}{2^{j/t}} \right]\right)  \\
&=& \Omega\left(m(\log d_{av})^{1/t}
    \left[\epsilon 2^{-\epsilon n/t}n^{1+\epsilon /t} +
    2^{-\epsilon n/t}n^{\epsilon/t}\frac{(1-2^{-(n-l-1)/t})}{1-2^{-1/t}}
      \right]\right) \\
&=& \Omega\left(m(\log d_{av})^{1/t}\times 2^{-\epsilon n/t}
    \left[\epsilon n^{1+\epsilon /t} +
    n^{\epsilon/t} \frac{(1-2^{-(n-l-1)/t})}{1-2^{-1/t}}
    \right]\right) \\
&=& \Omega\left(m(\log d_{av})^{1/t}\times 2^{-\epsilon n/t}
    (\epsilon n^{1+\epsilon /t}+ \Theta(n^{\epsilon/t}))\right) \;\;\ldots\;\; (B)
\end{eqnarray*}
The ratio of the bound in $(B)$ to the one in $(A)$ can be seen to be
$\Omega(\left(2^n/n\right)^{(1-\epsilon)/t})$, which is
$\Omega((\log |V|/\log \log |V|))^{(1-\epsilon)/t})$
for an appropriately slow-growing function $w$.

\end{enumerate}

\section{Concluding Remarks} \label{sec:concln}

In the course of this paper, some semi-combinatorial
proofs of certain non-trivial identities involving
binomial coefficients were also obtained. These are
described below:

\begin{eqnarray}
\sum_{a=0}^A \sum_{i=0}^{d-a} {d\choose {a+i}}{a+i\choose i}2^i(2d-2a-i)!(2a+i)! &=& (d!)^24^{d-A}(A+1){2A+1 \choose A}
\end{eqnarray}

The LHS (when divided by $(2d+1)!$) amounts to the expression for
$P_v(A)$ when $k=3$: choose $a+i$ hyperedges from the $d$ hyperedges
incident on $v$, of these $a$ hyperedges are backward,
while $i$ hyperedges each have one vertex occurring
prior to $v$ in the random permutation. These $i$
vertices can be chosen from $i$ pairs in $2^i$ ways.
The $(2a+i)$ vertices before $v$ can be arranged
in $(2a+i)!$ ways amongst themselves.
The remaining $(2d-2a-i)$ vertices occur after $v$
and can be arranged amongst themselves in $(2d-2a-i)!$
ways. The RHS is easily obtained from Theorem
\ref{genrlathm} by taking $t=2$.

\paragraph{}
Taking $A=0$ in the above expression gives us the simpler identity:

\begin{eqnarray}
\sum_{i=0}^d {d\choose i}2^i(2d-i)!i! & = & (d!2^d)^2
\end{eqnarray}

Dividing by $(d!)^2 2^d$ and changing the order of
summation, we get

\begin{eqnarray}
\sum_{i=0}^d {d+i \choose d}2^{-i} &=& 2^d
\end{eqnarray}

The above expression is discussed in some detail in
\cite{GrKnPa94} (Chapter 5, eqs. 5.20, 5.135-8); a nice combinatorial
proof of it is provided in \cite{Tamas93}.

\paragraph{}
The next expression (for the more general case $k \geq 3$) is
much more complicated:

\begin{eqnarray}
\sum_{a=0}^A \sum_{i=0}^{d-a} \sum_{\sum_{j=1}^{t-1}i_j=i;i_j\geq 0}
{d\choose {a+i}}{a+i\choose {a,i_1,\ldots,i_{t-1}}}
{t\choose 1}^{i_1}{t\choose 2}^{i_2}\ldots
{t\choose t-1}^{i_{t-1}}\times \nonumber \\
(ta+i_1+2i_2+\ldots (t-1)i_{t-1})!(td-ta-i_1-\ldots-(t-1)i_{t-1})!  \nonumber \\
 = (td+1)!(1+(tA+t)^{-1})^{-1}\frac{{d\choose A+1}}{{d+1/t\choose d-A-1}}
\end{eqnarray}

The LHS again follows by similar arguments as in the
previous case, this time for general $t$. There are
$a$ backward edges, $i_1$ edges which
have one vertex before $v$, $i_2$ edges with 2
vertices before $v$, and so on. The RHS follows
from Theorem \ref{genrlathm}.

\vspace{2mm}
It was in fact the non-triviality of the above LHS
expressions that led to our use of the PIE and its
variant (Theorem \ref{mpiethm}) in order to obtain
closed-form expressions for Theorems \ref{azerothm}
and \ref{genrlathm}.

\paragraph{}
With regard to the tightness of our results and the
weakening parameter $A$, firstly, from the proof of
Theorems \ref{trianglefreethm}-\ref{kunithm}, it is
clear that $\epsilon = \log A/\log D$ has to be at
least a constant. Ideally, we may want to replace $A$ by 1 in the
bounds of Theorems \ref{trianglefreethm}, \ref{krfreethm} and
\ref{kunithm}. This corresponds to the case $\epsilon = 0$. The
following example, however, shows that it is possible to construct
a triangle-free graph for which the bound in say, Theorem
\ref{trianglefreethm} would give a value more than the number
of vertices: Take a disjoint union of $A=K_{n/3,n/3}$ and
$B=\overline{K}_{n/3}$, and introduce a perfect matching between $B$ and
one of the parts of $A$. Now, $|V|= n$, $D\sim 2n/9$,
and hence if $A=1$, Theorem
\ref{trianglefreethm} would give a lower bound of
$\Omega(n\log n)$, which is asymptotically larger than
$|V|$. Similar examples can be constructed with linear
hypergraphs also.

\paragraph{{\bf Acknowledgement}} The first and third authors
would like to thank N.R. Aravind for some initial helpful
discussions.

\section{Appendix}

\begin{proof}[Proof of Lemma \ref{azerolemma}]\footnote{A proof based on $n$-th order difference operators follows from
\cite{GrKnPa94} (Chapter 5, eqn. 5.41). However, we were
not aware of this at the time of solving. The hypergeometric
proof is equally simple and the approach more general.}
Write the LHS as $\sum_{r\geq 0} t_r$, since ${d\choose r} = 0$ for
$r >d$. Now,
\begin{eqnarray*}
   \frac{t_{r+1}}{t_r} &=& \frac{(-1)(d-r)(tr+1)}{(r+1)(tr+t+1)} \\
                       &=& \frac{(r-d)(r+1/t)}{(r+1)(r+1+1/t)} \\
\end{eqnarray*}
Also, notice that $t_0=1$. Therefore, the LHS can be written as
the generalised hypergeometric function $F(1/t,-d;1+1/t;1)$, where the
generalised hypergeometric function $F(a_1,...,a_m;b_1,...,b_n;z)$ is
given by
\begin{eqnarray*}
  F(a_1,...,a_m;b_1,...,b_n;z) &=& \sum_{r=0}^{\infty} \frac{(a_1)^{(r)}(a_2)^{(r)}...(a_m)^{(r)}}{(b_1)^{(r)}...(b_n)^{(r)}}\frac{z^r}{r!}
\end{eqnarray*}
where $p^{(q)} = p(p+1)...(p+q-1)$ is the \emph{rising factorial}.
Next, we use the general version of Vandermonde convolution - also
known as Chu-Vandermonde identity (a special case of
Gauss's Hypergeometric Theorem, see e.g. \cite{GrKnPa94},
Chapter 5, equation 5.93, also \cite{AnAsRo99,Ba35,Be02,Wi})
$$F(a,-n;c;1) = \frac{(c-a)^{(n)}}{c^{(n)}}$$
The above is true whenever $a,c$ are complex numbers and $n$ is a natural
number, such that $\mathbb{R}(a)-n< \mathbb{R}(c)$.
In our case, $a = 1/t$, $n=d$ and $c= 1+1/t$. Hence we get $(c-a)^{(n)}=d!$,
and $c^{(n)} = (1+1/t)^{(d)} = (1+1/t)(2+1/t)...(d+1/t)$. Therefore, the
LHS of (\ref{binident1}) becomes
$$F(1/t,-d;1+1/t;1) = \frac{d!}{(1+1/t)(2+1/t)...(d+1/t)} =
\frac{1}{{{d+1/t}\choose d}}\qquad \qquad
\hfill \Box $$
\end{proof}

\begin{proof} [Proof of Lemma \ref{mpielemma}]
  The proof is by induction on $b$. For $b=0$, the LHS reduces to
$$\sum_{i=0}^0 (-1)^i {{a+0}\choose {a+i}}{{a+i-1}\choose i}$$ which is
clearly $1$.
Assume the lemma to be true for $b=c$ and consider the LHS when $b=c+1$:
\begin{eqnarray*}
& & \sum_{i=0}^{c+1} (-1)^i {{a+1+c}\choose {a+i}}{{a+i-1}\choose i} \\
&=& \sum_{i=0}^{c+1} (-1)^i \left[{{a+c}\choose {a+i}}+{{a+c}\choose {a+i-1}}
                            \right]{{a+i-1}\choose i} \\
&=& \sum_{i=0}^{c+1} (-1)^i \left[{{a+c}\choose {a+i}}{{a+i-1}\choose i}+
                            {{a+c}\choose {a+i-1}}{{a+i-1}\choose i}\right] \\
&=& 1+ \sum_{i=0}^{c+1} (-1)^i {{a+c}\choose {a+i-1}}{{a+i-1}\choose i}
\end{eqnarray*}
by the induction hypothesis, since ${{a+c}\choose {a+c+1}} = 0$.
Now, the second sum is
\begin{eqnarray*}
& & \sum_{i=0}^{c+1} (-1)^i{{a+c}\choose {a+i-1}}{{a+i-1}\choose i} \\
&=& \frac{(a+c)!}{(a-1)!(c+1)!}\sum_{i=0}^{c+1} (-1)^i{{c+1}\choose i} \\
&=& 0
\end{eqnarray*}
\hfill $\Box$
\end{proof}

\begin{proof}[Proof of Lemma \ref{genrlaexprlemma}]
Let the LHS be denoted by $S_d$. Then, using the identity
${n\choose r} = {{n-1}\choose r}+{{n-1}\choose{r-1}}$, we have
\begin{eqnarray*}
S_d &=& \sum_{r=0}^{d-A} (-1)^r \left[{{d-1}\choose{r+A}}+
  {{d-1}\choose{r+A-1}}\right]{{A+r-1}\choose r}\frac{1}{tr+tA+1}\\
    &=& S_{d-1} + \sum_{r=0}^{d-A} (-1)^r {{d-1}\choose{r+A-1}}
       {{A+r-1}\choose r}\frac{1}{tr+tA+1}\\
\end{eqnarray*}
since ${{d-1}\choose d} = 0$. Now the second sum can be simplified as
\begin{eqnarray*}
T_d &=& \sum_{r=0}^{d-A} (-1)^r {{d-1}\choose{r+A-1}}{{A+r-1}\choose r}
                              \frac{1}{tr+tA+1} \\
&=& \left(\frac{(d-1)!}{(d-A)!(A-1)!}\right)
    \sum_{r=0}^{d-A} (-1)^r {{d-A}\choose{r}}\frac{1}{tr+tA+1} \\
&=& {{d-1}\choose{A-1}}\frac{1}{tA+1}
    \sum_{r=0}^{d-A} (-1)^r {{d-A}\choose{r}}\frac{1}{(t/(tA+1))r+1} \\
\end{eqnarray*}
By Lemma \ref{azerolemma}, we get
\begin{eqnarray*}
T_d &=& \frac{1}{tA+1}\frac{{{d-1}\choose{A-1}}}{{{d-A+(tA+1)/t}\choose {d-A}}}
\end{eqnarray*}
Therefore,
\begin{eqnarray*}
S_d &=& S_{d-1} + \frac{1}{tA+1}\frac{{{d-1}\choose{A-1}}}{{{d+1/t}\choose {d-A}}}
\end{eqnarray*}
Unraveling the recursion and noticing that $S_A = 1/(tA+1)$, we get that
\begin{eqnarray*}
S_d &=& (1/(tA+1))\sum_{r=0}^{d-A} \frac{{{d-r-1}\choose{A-1}}}
                                        {{d+1/t-r \choose d-A-r}}\\
&=& (1/(tA+1))\sum_{r=0}^{d-A} \frac{{{A-1+r}\choose{A-1}}}
                                    {{A+1/t+r \choose r}}
\end{eqnarray*}
by reversing the order of summation.
Finally, the following claim completes the
proof.
\medskip

{\bf Claim.}
For $d\geq A, t\geq 0$,
$$\frac{1}{tA+1}\sum_{r=0}^{d-A} \frac{{{A-1+r}\choose{A-1}}}
  {{A+1/t+r \choose r}} \;\;=\;\; 1-\frac{tA}{tA+1}\frac{{d \choose A}}{{d+1/t \choose d-A}}$$

{\bf Proof of Claim.} We use induction on $d$. When $d=A$, the LHS is
$(tA+1)^{-1}$, while the RHS is $1-\frac{At}{tA+1}$, so we have equality.
Now assume equality for $d$ and consider the LHS for $d+1$:
\begin{eqnarray*}
& & \frac{1}{tA+1}\sum_{r=0}^{d-A+1}\frac{{A-1+r\choose r}}{{A+1/t+r\choose r}}
    \\
&=& 1-\frac{At}{tA+1}\frac{{d\choose A}}{{d+1/t\choose d-A}}+(At+1)^{-1}\frac{{d\choose d-A+1}}{{d+1+1/t\choose d-A+1}} \\
&=& 1-\frac{At}{(tA+1){d+1+1/t\choose d-A+1}}\left[{d\choose A}\frac{d+1+1/t}{d-A+1}-(At)^{-1}{d\choose d-A+1}\right]\\
&=& 1-\frac{At}{(tA+1){d+1+1/t\choose d-A+1}}\left[{d+1\choose A}+{d+1\choose d-A+1}(t(d+1))^{-1}-(At)^{-1}{d\choose A-1}\right]\\
&=& 1-\frac{At}{(tA+1){d+1+1/t\choose d-A+1}}{d+1\choose A}
\end{eqnarray*}
which is the required expression on the RHS.
\hfill $\Box$
\end{proof}

\end{document}